\documentclass[12pt]{article}

\usepackage{graphicx}
\usepackage{amsfonts}
\usepackage{amsthm}
\usepackage{euscript}

\usepackage{caption}

\usepackage[colorlinks, citecolor=black, urlcolor=blue, linkcolor=black]{hyperref}

\newtheorem{theorem}{Theorem}[section]

\makeatletter
\@addtoreset{equation}{section}
\makeatother

\usepackage{fancyhdr}
\begin{document}

\thispagestyle{plain} 

\sloppy

\begin{center}
{\bf \Large\bf Spontaneous decay of level from spectral theory point of view.}
\\ \bigskip
{\large Ianovich E.A.}
\\ \bigskip
   {\it Saint Petersburg, Russia}
\\ \bigskip
   {\it\small E-mail: eduard@yanovich.spb.ru}
\end{center}

\begin{abstract}
In quantum field theory it is believed that the spontaneous decay of excited atomic or molecular level is due to the interaction with continuum of field modes. Besides, the atom makes a transition from upper level to lower one so that the probability to find the atom in the excited state tends to zero. In this paper it will be shown that the mathematical model in single-photon approximation may predict another behavior of this probability generally. Namely, the probability to find the atom in the excited state may tend to a nonzero constant so that the atom is not in the pure state finally. This effect is due to that the spectrum of the complete hamiltonian is not purely absolutely continuous and has a discrete level outside the continuous part. Namely, we state that in the corresponding invariant subspace, determining the time evolution, the spectrum of the complete hamiltonian when the field is considered in three dimensions may be not purely absolutely continuous and may have an eigenvalue. The appearance of eigenvalue has a threshold character. If the field is considered in two dimensions the spectrum always has an eigenvalue and the decay is absent.
\end{abstract}

\bigskip

\section{Introduction}

The formal physical energy operator (hamiltonian) of the two-level system ("atom") interacting with quantized field in single-photon approximation has the form
\begin{equation}
\label{hamiltonian}
H=E_2\,|2\rangle\langle 2|+E_1\,|1\rangle
\langle1|+\int\limits_{0}^{+\infty}x\,a^+_x a_x\,dx+\int\limits_{0}^{+\infty}\left\{V(x)
|2\rangle\langle1|\,a_x+V^*(x)|1\rangle\langle2|\,a^+_x\right\}dx\,,
\end{equation}
where
$$
H_0=E_2\,|2\rangle\langle 2|+E_1\,|1\rangle\langle1|+\int\limits_0^{+\infty}x\,a^+_x a_x\,dx
$$
is unperturbed hamiltonian of free two-level system and field and
$$
H_{int}=\int\limits_0^{+\infty}\left\{V(x)
|2\rangle\langle1|\,a_x+V^*(x)|1\rangle\langle2|\,a^+_x\right\}dx
$$
is the hamiltonian of interaction so that $H=H_0+H_{int}$. Here we used Dirac's notations: $|1\rangle$ and $|2\rangle$ are states of two-level system with energies $E_2>E_1$. Operators $a^+_x$ and $a_x$ are operators of creation and annihilation of photon with energy $x$, satisfying the standard commutation relations
$$
[a_x,a^+_{x'}]=\delta(x-x')
$$

We neglect polarization of photons and consider quantized field as a scalar boson field in one dimension. In fact, such a model is sufficient for consideration the spontaneous decay as well as the spectrum of hamiltonian in two or real three-dimensional case. Namely, one, two and three-dimensional cases correspond to a special choice of the function $V(x)$. In this work we consider two cases: three-dimensional and two-dimensional. In both cases the function $V(x)$ is supposed to be continuous on $[0;+\infty)$ and have no zeros in the open interval $(0;+\infty)$. Moreover, in both cases $V(x)\in L_2[0;\infty)$ so that
$$
\int\limits_0^{+\infty}|V(x)|^2dx<\infty
$$

These conditions are fulfilled in the most physical applications. The condition $|V(x)|>0$ for $x>0$ is not prohibited by any physical law but converse is rather exception than a rule.

The principal difference between these two cases is concluded in the behaviour of function $V(x)$ in the point $x=0$.
In three-dimensional case $V(0)=0$ so that $|V(x)|^2=xV_1(x)$, where $V_1(0)\ne0$~\cite{4}. In two-dimensional case $V(0)\ne0$.
This difference is due to the difference in jacobian in spherical and polar coordinates.
Note also that in one-dimensional case $V(x)\notin L_2[0;\infty)$.

n-photon states of field are defined by
$$
|x_1,\ldots,x_n\rangle=
a^+_{x_1}\ldots a^+_{x_n}\,|vac\rangle\,,
$$
where $|vac\rangle$ is a vacuum state of field: $a_{x}|vac\rangle=0$.

Suppose at time $t=0$ the whole system is in the state $|\Psi(0)\rangle=|2\rangle|vac\rangle$. That is the atom is in the excited state and the field is in the vacuum state. In single-photon approximation the probability $P(t)$ to find the atom in the excited state at moment $t$ is defined by the quantum amplitude
\begin{equation}
\label{quantum_amplitude}
C(t)=\langle\Psi(0)|\,e^{-iHt}\,|\Psi(0)\rangle
\end{equation}
so that
\begin{equation}
\label{probability}
P(t)=|C(t)|^2.
\end{equation}
Note that this amplitude satisfies the following integro-differential equation~\cite{1}
$$
\frac{dy(t)}{dt}=\int\limits_0^t K(t-t')\,y(t')\,dt'\,,
$$
where
$$
y(t)=C(t)e^{iE_2t}\,,\quad K(t)=-e^{i(E_2-E_1)t}\int\limits_0^{+\infty}|V(x)|^2\,e^{-ixt}\,dx.
$$

As it is believed~\cite{1,2} the absolute value of amplitude $|C(t)|$ tends to zero almost exponentially. The analysis was based on the consideration of this integro-differential equation.

We will apply here another method based on spectral properties of the hamiltonian $H$. We will show that in three-dimensional case the probability $P(t)$ tends to zero for relatively small interaction only. In quantum electrodynamics the role of this smallness parameter belongs to thin structure constant $\alpha=e^2/\hbar c\sim1/137$. For strong interaction as well as for two dimensional case at any strength of interaction the probability $P(t)$ tends to a nonzero constant when $t\to+\infty$, that is
$$
\lim\limits_{t\to+\infty}P(t)=P_\infty>0.
$$

As it will be proved below this effect is due to that the spectrum of the energy operator is not purely absolutely continuous in the invariant subspace defining this time evolution. In the two-dimensional case there exists an eigenvalue for any strength of interaction. In the three dimensional case the appearance of eigenvalue has threshold character.

\section{Spectrum of the hamiltonian in the corresponding invariant subspace}

It is easy to see from~(\ref{hamiltonian}) that the states $|2\rangle|vac\rangle$, $|1\rangle|x\rangle$ ($0\le x<+\infty$) form invariant subspace for the operator $H$. Mathematically this subspace is equivalent to the Hilbert space $\EuScript{H}=\EuScript{H}_1\bigoplus\EuScript{H}_2$, where $\EuScript{H}_1$ is one-dimensional space connected with discrete state $|2\rangle|vac\rangle$, and $\EuScript{H}_2$ is usual $L_2[0,+\infty)$ space connected with continuum of atomic-field states $|1\rangle|x\rangle$. Any vector $f\in\EuScript{H}$ can be presented in the form
$$
f=\left(
    \begin{array}{c}
      f_1 \\
      f_0
    \end{array}
  \right)\,,
$$
where $f_0\in\EuScript{H}_1$ ($f_0$ is constant), $f_1\equiv f_1(x)\in\EuScript{H}_2$. Scalar product in $\EuScript{H}$ is defined by
$$
(f,g)=f_0g^*_0+\int\limits_0^{+\infty} f_1(x)g^*_1(x)\,dx.
$$
Let $A$ be the operator generated by the operator $H$ in $\EuScript{H}$. Then the action of $A$ on the vector $f$ can be written as follows
\begin{equation}
\label{operator_A}
Af=A\left(
    \begin{array}{c}
      f_1 \\
      f_0
    \end{array}
  \right)=\left(
    \begin{array}{c}
      (x+E_1)f_1(x)+V^*(x)f_0 \\
      E_2f_0+\int\limits_0^{+\infty} V(x)f_1(x)\,dx
    \end{array}
  \right)
\end{equation}

The operator $A$ is defined on the domain $D_A=\{f\in H\,|\,(xf_1(x))\in L_2[0,\infty)\}$. It is easy to check that $A$ is self-adjoint. Note that the action of the operator $A_0$ generated by the operator $H_0$ in $\EuScript{H}$ is defined by
\begin{equation}
\label{A_0f}
A_0f=\left(
    \begin{array}{c}
      (x+E_1)f_1(x)\\
      E_2f_0
    \end{array}
  \right).
\end{equation}

It follows that the operator $A_0$ in the space $\EuScript{H}_2$ is the multiplication operator and has purely absolutely continuous spectrum on semi-axis $[E_1,+\infty)$. In the space $\EuScript{H}_1$ the operator $A_0$ spectrum consists of eigenvalue $E_2$ only.

To investigate the spectrum of $A$ one can use the resolvent identity
\begin{equation}
\label{resolvent_identity}
R_0(\lambda)=R(\lambda)+R(\lambda)VR_0(\lambda)\,,
\end{equation}
where $R_0(\lambda)=(A_0-\lambda E)^{-1}$, $R(\lambda)=(A-\lambda E)^{-1}$, $V=A-A_0$ and
\begin{equation}
\label{Vf}
    Vf=\left(
    \begin{array}{c}
      V^*(x)f_0\\
      \int\limits_0^{+\infty} V(x)f_1(x)\,dx
    \end{array}
  \right).
\end{equation}
Let
$$
f^{(1)}=\left(
    \begin{array}{c}
      0\\
      1
    \end{array}
  \right)\,,\quad
  f^{(2)}=\left(
    \begin{array}{c}
      V^*(x)\\
      0
    \end{array}
  \right)\,,\quad R_{ij}(\lambda)=(R(\lambda)f^{(i)},f^{(j)})\,,\:\:i,j=1,2.
$$
From resolvent identity~(\ref{resolvent_identity}) using~(\ref{A_0f}) and~(\ref{Vf}) one can find
$$
\left\{\begin{array}{c}
  R_{12}(\lambda)+\frac{1}{E_2-\lambda}R_{22}(\lambda)=0 \\
  k(\lambda)R_{12}(\lambda)+R_{22}(\lambda)=k(\lambda)
\end{array}\right.\,,
$$
$$
\left\{\begin{array}{c}
  k(\lambda)R_{11}(\lambda)+R_{21}(\lambda)=0 \\
  R_{11}(\lambda)+\frac{1}{E_2-\lambda}R_{21}(\lambda)=\frac{1}{E_2-\lambda}
\end{array}\right.\,,
$$
where
$$
k(\lambda)=\int\limits_0^{+\infty}\frac{|V(x)|^2\,dx}{x+E_1-\lambda}\,,
$$
whence
\begin{equation}
\label{R_{11}}
R_{11}(\lambda)=\frac{1}{E_2-\lambda-k(\lambda)}\,,\quad R_{22}(\lambda)=(E_2-\lambda)k(\lambda)\,R_{11}(\lambda)
\end{equation}

From~(\ref{R_{11}}) it follows that the functions $R_{11}(\lambda)$ and $R_{22}(\lambda)$ describe the same spectrum. We will see below that the spectrum of $A$ is simple and the vector $f^{(1)}$ is generating vector. The vector $f^{(2)}$ is also generating vector if $f^{(2)}\in D_A$. So, we can consider the function $R_{11}(\lambda)$ only. Let $t\in\mathbb{R}$. From~(\ref{R_{11}}), using Sokhotski-Plemelj theorem, we obtain
\begin{equation}
\label{limit_R_{11}}
\lim\limits_{\epsilon\to+0}R_{11}(t+i\epsilon)=\frac{1}{E_2-t-i\pi|V(t-E_1)|^2-
\displaystyle\EuScript{P}\int\limits_0^{+\infty}\frac{|V(x)|^2\,dx}{x+E_1-t}}
\end{equation}

Here $V(t-E_1)\equiv0$, when $t<E_1$. The integral in the denominator is understood in the sense of principal value in the point $t$ if $t\in(E_1;+\infty)$ and as usual improper integral if $t\le E_1$.
From~(\ref{limit_R_{11}}) it follows that if $V(t-E_1)\ne 0$ at some point $t$ then this point together with some their neighbourhood belongs to the absolutely continuous part of the operator $A$ spectrum (remember that $V(x)$ is continuous). Moreover, if $\sigma(t)$ is the distribution function connected with integral representation of $R_{11}(\lambda)$
$$
R_{11}(\lambda)=\int\limits_{-\infty}^{+\infty}\frac{d\sigma(t)}{t-\lambda}\,,
$$
then almost everywhere
$$
\sigma'(t)=\sigma_{ac}'(t)=\rho(t)=\frac{1}{\pi}\,\mbox{Im}\left[\lim\limits_{\epsilon\to+0}R_{11}(t+i\epsilon)\right]=
$$
\begin{equation}
\label{density}
\frac{|V(t-E_1)|^2}
{\left(E_2-t-\displaystyle\EuScript{P}\int\limits_0^{+\infty}\frac{|V(x)|^2\,dx}{x+E_1-t}\right)^2+\pi^2|V(t-E_1)|^4}
\end{equation}
Here $\sigma_{ac}(t)$ is the absolutely continuous part of $\sigma(t)$.

If the point $t$ belongs to the singular spectrum of $A$, then
$$
V(t-E_1)=0\,,\quad E_2-t-\displaystyle\EuScript{P}\int\limits_0^{+\infty}\frac{|V(x)|^2\,dx}{x+E_1-t}=0
$$
Note that a very close situation takes place in Friedrichs model~\cite{3}.

\begin{figure}[h]
  \center{\includegraphics[scale=0.5]{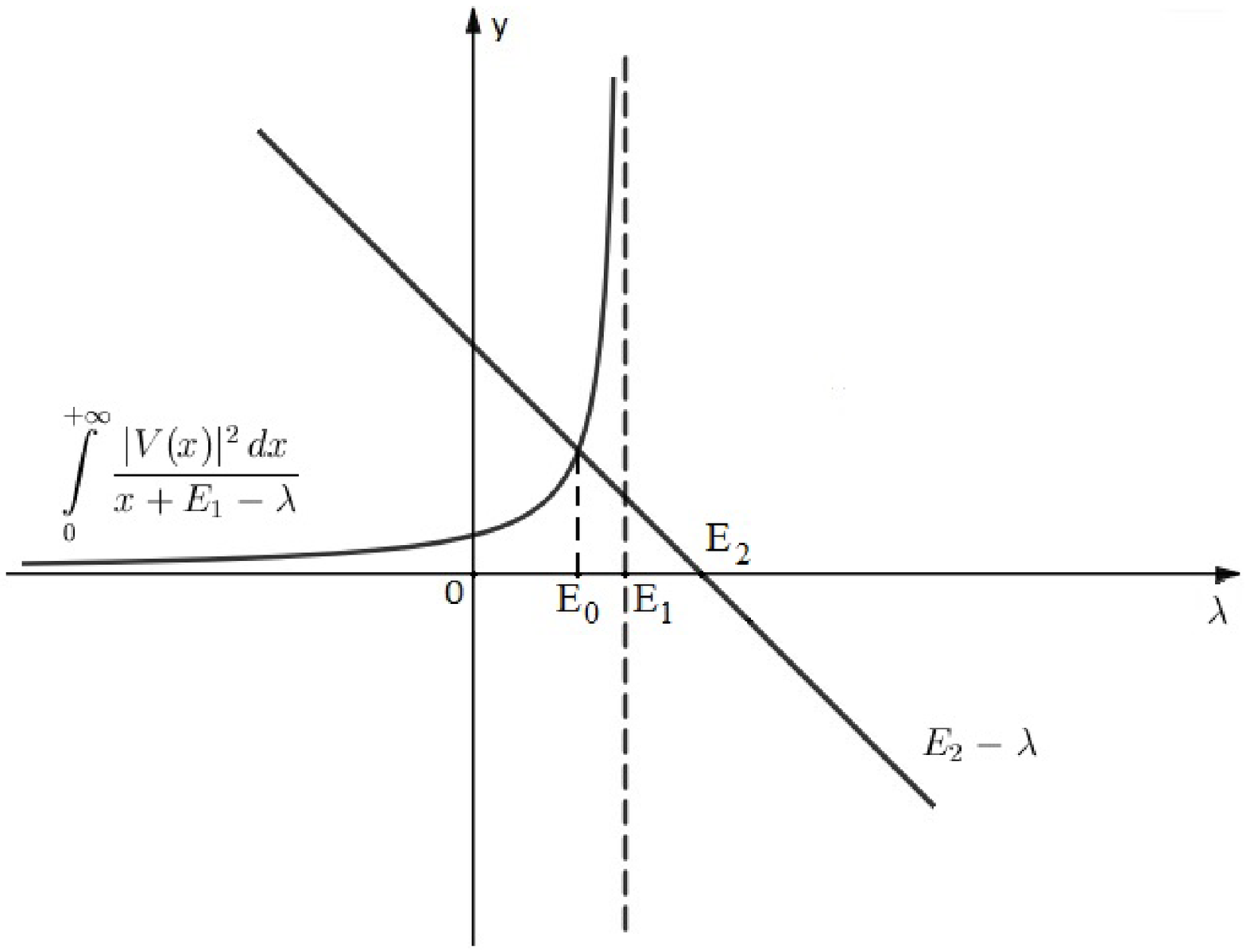}}
  \caption{}
  \label{Fig.1}
\end{figure}

Since we suppose $|V(x)|>0$ for $x>0$, the absolutely continuous spectrum of $A$ covers semi-axis $[E_1;+\infty)$. A singular part of spectrum may be not empty for $\lambda\le E_1$ only. Let us consider first the two-dimensional case. Since $V(0)\ne 0$, the singular spectrum may exist for $\lambda<E_1$ only and consists of isolated eigenvalues which are roots of the equation
\begin{equation}
\label{eigen_equation}
E_2-\lambda=\int\limits_0^{+\infty}\frac{|V(x)|^2\,dx}{x+E_1-\lambda}\quad (\lambda<E_1)
\end{equation}
The Figure~(\ref{Fig.1}) shows a sketch of the graphs of the right and left hand sides of this equation.

From this figure it follows that the equation~(\ref{eigen_equation}) has a unique solution $E_0$ such that $E_0<E_1$. $E_0$ is an eigenvalue of the operator $A$. Actually, the equation $Af=\lambda f$ reduces to the system
$$
\left\{\begin{array}{c}
         (x+E_1)f_1(x)+V^*(x)f_0=\lambda f_1(x) \\
          E_2 f_0+\int\limits_0^{+\infty}V(x)f_1(x)\,dx=\lambda f_0
       \end{array}\right.
$$

Expressing $f_1(x)$ from the first equation and substituting it to the second one, we obtain the following equation for $\lambda$
$$
E_2+\int\limits_0^{+\infty}\frac{|V(x)|^2\,dx}{\lambda-x-E_1}=\lambda\,,
$$
which coincides with~(\ref{eigen_equation}). Note that since for any $\lambda$ the last system has exactly one solution with accuracy of a multiplication constant $f_0$, the spectrum of the operator $A$ is simple. The eigenvector $e$, corresponding to the eigenvalue $E_0$, has the form
\begin{equation}
\label{eigenvector}
e=\left(\begin{array}{c}
            \displaystyle -\frac{V^*(x)}{E_1+x-E_0} \\
            1
          \end{array}
\right)\,,
\end{equation}
and
\begin{equation}
\label{norm_of_eigenvector}
\|e\|^2=1+\int\limits_0^{+\infty}\frac{|V(x)|^2\,dx}{(x+E_1-E_0)^2}<+\infty\,,
\end{equation}
since $E_0<E_1$. Hence $e\in D_A$.

Let us consider now the three-dimensional case. Since $|V(x)|^2=xV_1(x)$ ($V_1(0)\ne0$), the integral in the right hand side of the equation~(\ref{eigen_equation}) converges at $\lambda=E_1$. From~(\ref{eigen_equation}) it follows that the eigenvalue at the
point $E_0<E_1$ exists if and only if
\begin{equation}
\label{porog}
E_2-E_1<\int\limits_0^{+\infty} V_1(x)\,dx=\int\limits_0^{+\infty}\frac{|V(x)|^2}{x}\,dx
\end{equation}

The existence of eigenvalue at the point $\lambda=E_1$ is impossible because in this case $\|e\|=+\infty$ and $e\notin\EuScript{H}$.

Resuming the results of this section, we obtain
\begin{theorem}
\label{spectrum_of_A}
Let $A$ be the operator~(\ref{operator_A}) induced by the hamiltonian~(\ref{hamiltonian}) in the corresponding invariant subspace. The operator $A$ is self-adjoint and its spectrum is simple. Suppose the function $V(x)$ is continuous on $[0;+\infty)$ and have no zeros in the open interval $(0;+\infty)$. Moreover, suppose $V(x)\in L_2[0;\infty)$, and one of the following two conditions\\

$1.\,V(0)\ne0$ (two-dimensional case)\\
or

$2.\,|V(x)|^2=xV_1(x)\,,\:V_1(0)\ne0$ (three-dimensional case)\\
is satisfied.

Then in the first case the operator $A$ has absolutely continuous spectrum on $[E_1;+\infty)$ and single eigenvalue $E_0<E_1$.
In the second case the spectrum of the operator $A$ is also absolutely continuous on $[E_1;+\infty)$ but the eigenvalue
$E_0<E_1$ appears if and only if the condition~(\ref{porog}) is fulfilled.
\end{theorem}

\section{Spontaneous decay low}

Knowing the spectrum of the operator $A$, we can analyse the behaviour of the quantum amplitude~(\ref{quantum_amplitude}) and hence the probability~(\ref{probability}). Let $E_{\lambda}$ be the resolution of the identity for the operator $A$. We have
$$
C(t)=\langle\Psi(0)|\,e^{-iHt}\,|\Psi(0)\rangle=(e^{-iAt}f^{(1)},f^{(1)})=
\int\limits_{-\infty}^{+\infty}e^{-i\lambda t}\,d(E_{\lambda}f^{(1)},f^{(1)})=\int\limits_{-\infty}^{+\infty}e^{-i\lambda t}\,d\sigma(\lambda)
$$
Using the results of the previous section, we obtain
$$
C(t)=\frac{|(e,f^{(1)})|^2}{\|e\|^2}\,e^{-iE_0 t}+\int\limits_{E_1}^{+\infty}e^{-i\lambda t}\,\rho(\lambda)\,d\lambda\,,
$$
where the density function $\rho(t)$ is defined by~(\ref{density}). Here we supposed that the spectrum contains eigenvalue $E_0$.
If it is not, then the first term is absent. By Rimann-Lebesgue lemma, the integral, corresponding to the absolutely continuous spectrum tends to zero as $t\to\infty$. So, we have
$$
C(t)=\frac{|(e,f^{(1)})|^2}{\|e\|^2}\,e^{-iE_0 t}+\circ(1)\,,\quad t\to\infty
$$
and
\begin{equation}
\label{limit_constant}
\lim\limits_{t\to+\infty}P(t)=\left(1+\int\limits_0^{+\infty}\frac{|V(x)|^2\,dx}{(x+E_1-E_0)^2}\right)^{-2}>0\,,
\end{equation}
where we used~(\ref{eigenvector}) and~(\ref{norm_of_eigenvector}).

Taking into account the Theorem~(\ref{spectrum_of_A}), we obtain
\begin{theorem}
\label{limit_of_probability}
In the two dimensional case and in the the three-dimensional case, when the condition~(\ref{porog}) is fulfilled, the probability $P(t)$ tends to a non-zero constant~(\ref{limit_constant}) as $t\to+\infty$. So, the decay is absent.

If in the three-dimensional case the condition~(\ref{porog}) is not fulfilled, then $P(t)\to 0$ as $t\to+\infty$.
\end{theorem}

\section{Conclusion}

It may to seem that the threshold condition~(\ref{porog}) always depends on the difference between atomic levels $(E_2-E_1)$. However in quantum electrodynamics it is not true. The function $V_1(x)$ is proportional to $(E_2-E_1)\,\alpha^2$, where $\alpha=\frac{\textstyle e^2}{\textstyle \hbar c}$ is thin structure constant~\cite{4}, so that the integral in the right hand side of~(\ref{porog}) is of the same order. Hence the threshold condition~(\ref{porog}) does not depend on $(E_2-E_1)$ and is not fulfilled, since $\alpha\sim1/137$. Thus in quantum electrodynamics the spontaneous decay of level is due to the smallness of the relative coupling constant $\alpha$.

Note also that the discrete part of the operator $H_0$ spectrum may undergo a jump at the perturbation $H_{int}$. Actually, the operator $A_0$ has an eigenvalue $E_2$, although the eigenvalue of the perturbed operator $A$ is $E_0$ and $E_0<E_1<E_2$ for any strength of perturbation $V$ in the two-dimensional case and for perturbation, satisfying threshold condition~(\ref{porog}), in the three-dimensional case.

Moreover, we can see that in the real three-dimensional case in quantum electrodynamics the spectrum of $A$ is purely absolutely continuous and the using of perturbation theory (which is used usually in physics~\cite{5,6,7}) has no any mathematical sense. So, such famous problems as Lamb's shift, mass renormalization and regularization, can be connected with fact that the complete hamiltonian has no eigenvalues at all.

\end{document}